\documentclass[a4paper,11pt,english]{elsarticle} 
\usepackage[utf8x,utf8]{inputenc}
\usepackage{lineno,hyperref}
\usepackage{color}
\hypersetup{colorlinks = true, linkcolor = red, citecolor = blue, bookmarksopen = true}
\usepackage{caption}

\usepackage{amsmath, dsfont} 
\usepackage{amsthm}
\usepackage[left=3cm,right=3cm,top=3cm,bottom=3.5cm]{geometry}
\usepackage{stmaryrd}
\usepackage{bbold} 
\usepackage{bm}
\usepackage{epstopdf}

\usepackage{booktabs}
\usepackage{array}

\usepackage{graphicx}
\usepackage[tight,footnotesize]{subfigure}

\usepackage[dvipsnames,svgnames,x11names]{xcolor}
\usepackage[ruled,vlined]{algorithm2e}

\usepackage{calc}
\usepackage{ulem}

\renewcommand\d{\textrm{d}}
\newcommand\CDM{\text{\textnormal{CDM}}}

\newcommand\X{\textbf{X}}

\newcommand\R{\mathds{R}}
\newcommand\p{\mathds{P}}
\newcommand\x{\textbf{x}}

\newcommand\TV{\textrm{TV}}

\newcommand\E{\mathds{E}}

\renewcommand{\d}{{\rm{d}}}
\newcommand{\supp}{{\rm{Supp}}}

\begin{document}

\begin{frontmatter}

\title{Simultaneous estimation of complementary moment independent and reliability-oriented sensitivity measures}

\author[a,b]{Pierre Derennes}
\ead{pierre.derennes@onera.fr}

\author[b]{J\'er\^ome Morio\corref{cor1}}
\ead{Jerome.Morio@onera.fr}
\cortext[cor1]{Corresponding author}

\author[c]{Florian Simatos}
\ead{florian.simatos@isae.fr}

\address[a]{Universit\'e de Toulouse, UPS IMT, F-31062 Toulouse Cedex 9, France}
\address[b]{ONERA/DTIS, Universit\'e de Toulouse, F-31055 Toulouse, France}
\address[c]{ISAE-SUPAERO and Universit\'e de Toulouse, Toulouse, France}

\begin{abstract}
In rare event analysis, the estimation of the failure probability is a crucial objective. However, focusing only on the occurrence of the failure event may be insufficient to entirely characterize the reliability of the considered system. This paper provides a common estimation scheme of two complementary moment independent sensitivity measures, allowing to improve the understanding of the system's rare event. Numerical applications are performed in order to show the effectiveness of the proposed estimation procedure.
\end{abstract}


\end{frontmatter}

\section{Introduction}

In diverse disciplines, systems modeling is often achieved by considering a black-box model for which the observation is expressed as a deterministic function of external parameters representing some physical variables. These basic variables are usually assumed random in order to take phenomenological uncertainties into account. Then, sensitivity analysis (SA) techniques play a crucial role in the handling of these uncertainties and in the comprehension of the system behavior. These techniques aim at identifying and ranking inputs with respect to their impact on the output. In addition, SA methods present two main objectives: decrease the output uncertainty by reducing uncertainty of the most influential inputs, and simplify the model by omitting contribution of least ones. The influence criterion depends on the considered SA approach. There are various SA techniques in literature and essentially two families stand out: local and global sensitivity analysis (GSA) methods, see \cite{iooss2015review,wei2015variable,borgonovo2016sensitivity,borgonovo2016common} and associated references for a review. Local methods aim at studying the behavior of the output locally around a nominal value of inputs. In contrast, global methods consider the whole variation range of inputs. 

Sensitivity analysis may also be performed with a rare event perspective. Reliability-oriented sensitivity analysis (ROSA) differs from the classical one in the nature of the output quantity of interest under study. Indeed, sensitivity analysis focuses on the model output whereas ROSA is broadly concerned with a reliability measure, typically the failure probability associated to an unsafe and undesired state of the system. Various ROSA methods have been proposed. First, several global sensitivity methods have been developed: for instance, a failure probability-based method \cite{cui2010moment}, variance decomposition-based methods \cite{wei2012efficient,yun2016efficient,yun2018efficient}, a method based on density perturbation \cite{lemaitre2015density} or, more recently, a quantile-oriented sensitivity approach \cite{browne2017estimate}. Local sensitivity method (often based on partial derivatives of the failure probability with respect to distribution parameters) are also available, see for instance \cite{vincent2018} for a comprehensive review of these methods. 

Recently, this scope has been expanded by \cite{raguet2018target} which proposes to classify ROSA methods in two different families:

\begin{itemize}
\item[$\bullet$] First, target (or regional) sensitivity analysis, which aims at studying the impact of inputs over a function of the output, typically the indicator function of a critical domain. 

\item[$\bullet$] Second, conditional sensitivity analysis, which aims at studying the impact of inputs exclusively within the critical domain, namely, conditionally to the failure event. 
\end{itemize}

These two points of view can lead to widely different answers. To illustrate this aspect, consider for instance the following simple toy model:
\begin{equation} \label{toy_model}
	Y = X_1 + \mathds{1}_{X_1>3} \lvert X_2 \rvert
\end{equation}
where $X_1$ and $X_2$ are independent centered Gaussian random variables with respective variance $1$ and $5$. The random variables $X_1$ and $X_2$ are viewed as the system input and $Y$ as the system output. Let us consider that for this system, $\{Y > 3\}$ is the failure event, and try to answer the following question: which out of $X_1$ and $X_2$ is more important from a rare event perspective? Actually, the answer depends on the viewpoint considered:
\begin{itemize}
	\item if one is interested in the \textbf{impact of the input on the failure occurring or not}, then of course $X_1$ is highly influential and $X_2$, that only kicks in $Y$ upon failure, plays no role;
	\item if one is now interested in the most influential input \textbf{upon failure} occurring, then $X_2$ should intuitively be more important than $X_1$ because of its higher variance.
\end{itemize}

In this paper we focus on moment-independent indices, that have recently attracted increasing attention in order to alleviate some of the limitations of classical variance-based indices. Our main message is that several of these indices can be efficiently estimated simultaneously with failure samples that can be generated, for instance, with one run of a sequential Monte Carlo or importance sampling. 

The rest of this paper is organized as follows. The section \ref{def_target_cond} aims at introducing two different ROSA indices, $\bar \eta_i$ and $\delta_i^f$, which are intrinsically linked to to GSA method of Borgonovo \cite{borgonovo2007new}. In Section~\ref{sec:estimation} we present our simultaneous estimation scheme for $\bar \eta_i$ and $\delta_i^f$ measures and numerical applications are performed in Section~\ref{appl} to assess its efficiency. Section~\ref{sub:generalization} discusses how this scheme can be extended to a more general context. Our estimation scheme relies on the maximum entropy method which is recalled in~\ref{ME principle}.

\section{Two complementary moment independent sensitivity measures \label{def_target_cond}}

In this paper we focus on Borgonovo's indices originally proposed in~\cite{borgonovo2007new}, although our method can be generalized to more general indices as discussed in Section~\ref{sub:generalization}. Let in the sequel $Z \mid Z'$ denote a random variable with random distribution the distribution of $Z$ conditioned on $Z'$ and a deterministic scalar function $\mathcal{M}: \mathds{R}^d \to \mathds{R}$. To measure the sensitivity of the output $Y = \mathcal{M}(\X)$ with respect to one of its input $X_i$, where $\X = (X_1, \ldots, X_d)$, Borgonovo~\cite{borgonovo2007new} proposed in the case where $(X_i, Y)$ is absolutely continuous with respect to Lebesgue measure the index
\begin{equation} \label{eq:orig-delta}
	\delta_i = \frac{1}{2} \E \left[ \left \lVert f_Y - f_{Y \mid X_i} \right \rVert_{L^1(\R)} \right],
\end{equation}
i.e., half the average of the $L^1$ distance between the density of $Y$ and the random density of $Y$ conditioned on $X_i$. If $X_i$ has a high influence on $Y$, the conditional density should be different from the non-conditioned one and $\delta_i$ should thus take large values. For further references and more details on $\delta$-sensitivity measures the reader can consult \cite{borgonovo2016sensitivity}.

In this paper we will adopt a more general definition of Borgonovo's index, which will make it possible to consider cases where $(X_i, Y)$ is not absolutely continuous with respect to Lebesgue measure. The motivation stems from considering the influence of $X_i$ not only on $Y$ but on possibly discrete functions of $Y$ such as $\mathds{1}_{Y > S}$, which captures the influence of $X_i$ on the failure occurring or not.

For this generalization, we see Borgonovo's index as a measure of dependency between $X_i$ and $Y$. Namely, let $d_\TV(Z_1, Z_2)$ denote the total variation distance between the distributions of the random variables $Z_1$ and $Z_2$. When $Z_1$ and $Z_2$ are absolutely continuous with respect to Lebesgue measure, we have $d_\TV(Z_1, Z_2) = \frac{1}{2} \lVert f_{Z_1} - f_{Z_2} \rVert_{L^1(\R)}$ and so we adopt the following generalization of Borgonovo's index:
\begin{equation} \label{eq:delta}
	\delta_i = \E \left[ d_\TV \left( Y, Y \mid X_i \right) \right] = d_\TV \left( (X_i, Y), (X_i, Y') \right)
\end{equation}
such that $Y'$ and $Y$ are independent and identically distributed random variables. The second equality holds when $(X_i, Y)$ is absolutely continuous with respect to some product measure $\lambda(\d x) \otimes \mu(\d y)$ (typically, $(X_i, Y)$ is absolutely continuous with respect to Lebesgue measure, or $X_i$ is and $Y$ is a discrete random variable).

In a rare event context, we are interested in the impact of $X_i$ not only on $Y$ but also on the occurrence of some rare event which we write $\{Y > S\}$. This means that we are interested in the influence of $X_i$ on the random variable $\mathds{1}_{Y>S}$: the corresponding generalized Borgonovo's index is therefore given by
\begin{equation} \label{eq:def-eta}
	\eta_i = \E \left[ d_\TV(\mathds{1}_{Y>S}, \mathds{1}_{Y>S} \mid X_i) \right] = \mathds{E}\left[ \left \lvert \p(Y>S) - \p(Y>S \mid X_i) \right \rvert \right]
\end{equation}
which is actually twice the index proposed in Cui et al.~\cite{cui2010moment}. One of the drawback of this index is that it is unnormalized as it is upper bounded by twice the rare event probability $2 \p(Y > S)$. To obtain a $[0,1]$-valued index, we use the relation
\begin{equation} \label{eq:eta-bar-eta}
	\eta_i = 2 \p(Y > S) \times d_\TV(X_i, X_i \mid Y > S)
\end{equation}
observed in~\cite{wang2018new} and that can be derived using Bayes' Theorem, to propose the $[0,1]$-valued index
\begin{equation} \label{eq:def-bar-eta}
	\bar \eta_i = d_\TV(X_i, X_i \mid Y > S) = \frac{1}{2} \left \lVert f_{X_i} - f_{X_i \mid Y > S} \right \rVert_{L^1(\R)}.
\end{equation}

Complementary to this approach, we may also be interested in the influence of $X_i$ upon failure, which corresponds to considering $\delta_i$ but when all the random variables involved are conditioned upon the failure $Y > S$. Thus, this conditional index, denoted by $\delta^f_i$, is given by
\begin{equation} \label{eq:def-delta^f}
	\delta^f_i = \E \left[ d_\TV( Y \mid Y>S, Y \mid \{Y>S, X_i\}) \right].
\end{equation}
When $(X_i, Y)$ is absolutely continuous, this is a particular case of~\eqref{eq:delta} and so if we denote by $(\tilde X_i, \tilde Y)$ a random variable distributed as $(X_i, Y)$ conditioned on $Y> S$, then we have in this case
\begin{equation} \label{eq:delta-f}
	\delta^f_i = \frac{1}{2} \left \lVert f_{\tilde X_i, \tilde Y} - f_{\tilde X_i} f_{\tilde Y} \right \rVert_{L^1(\R^2)}.
\end{equation}

Instead of focusing on $Y$, the indices $\eta_i$ and $\bar \eta_i$ target a different output, namely $\mathds{1}_{Y > S}$ and will thus be referred as target indices. Similarly, instead of working in the normal mode, the indices $\delta^f_i$ are concerned with the system conditioned upon failure and will thus be referred to as conditional indices. See Section~\ref{sub:generalization} for more on this terminology.

For the toy model~\eqref{toy_model}, we have $Y > S$ if and only if $X_1 > S$: this directly implies $\p(Y > S \mid X_1) = \mathds{1}_{X_1 > S}$ and $\p(Y > S \mid X_2) = \p(Y > S)$ and then
\[ \bar \eta_1 = 1 - \p(X_1 > S) \approx 0.9987 \ \text{ and } \ \bar \eta_2 = 0. \]

This confirms the intuition that, as far as we are concerned with the failure occurring or not, $X_1$ is highly influential and $X_2$, not at all. However, in this simple Gaussian case we can directly compute the $\delta^f_i$'s through numerical integration, which gives
\[ \delta_1^f \approx 0.0781 \text{ and } \delta_2^f \approx 0.7686. \]

Thus upon failure, $X_2$ has become much more influential than $X_1$. This simple toy example illustrates the complementarity of the indices $\bar \eta_i$ and $\delta^f_i$ from a rare event perspective, and our goal in this paper is to show how they can be simultaneously and accurately estimated with only one run of sequential Monte Carlo or importance sampling, regularly considered in the context of rare event probability estimation $\p(Y > S)$. In other words, we show that upon estimating this probability, one also gets ``for free'', that is without additional calls to the function $\mathcal{M}$, an estimation of $\bar \eta_i$ and $\delta^f_i$.

\section{Simultaneous estimation of $\delta^f_i$ and $\bar \eta_i$ \label{sec:estimation}}

We consider throughout this article a general computer code $Y=\mathcal{M}(\X)$ where the scalar output $Y$ depends on a $d$-dimensional real valued random vector $\X=(X_1,\ldots,X_d)$ of $\mathds{R}^d$ through a deterministic scalar function $\mathcal{M}: \mathds{R}^d \to \mathds{R}$ called ``black box''. Without loss of generality, it is assumed that the failure event corresponds to the exceeding of a critical threshold $S$ by the output $Y$, i.e., is of the form $\{Y>S\}$.

We further assume that for every $i$, $(X_i, Y)$ is absolutely continuous with respect to Lebesgue measure with density $f_{X_i, Y}$ and marginals $f_{X_i}$ and $f_Y$. As above, we denote by $\tilde \X = (\tilde X_1, \ldots, \tilde X_d)$ a random variable distributed as $\X$ conditioned on $Y> S$ and define $\tilde Y = \mathcal{M}(\tilde \X)$. Thus, $(\tilde X_i, \tilde Y)$ is also absolutely continuous with respect to Lebesgue measure with density $f_{\tilde X_i, \tilde Y}$ with marginals $f_{\tilde X_i}$ and $f_{\tilde Y}$. Our simultaneous estimation scheme is obtained by combining state-of-the-art estimation techniques which we recall next.

\subsection{Estimation of $\delta_i$}

Initial estimations of $\delta$-sensitivity measures relied on their original definition in terms of total variation distance between conditional and unconditional distributions. Involving $L^1$ norms of differences of conditional and unconditional output probability density functions, this approach typically necessitates expensive double-loop estimation procedures with a prohibitive cost~\cite{borgonovo2007new, liu2009new, plischke2013global}. Alternative approaches were proposed in~\cite{zhang2014new, zhang2013structural}, but these two methods rest on strong technical assumptions such as independence between input or approximation of the black box $\mathcal{M}$ within the cut-HDMR (high-dimensional model representation) framework. An apparently efficient single-loop method was proposed in~\cite{wei2013monte}, but simulation results provided in~\cite{derennes2019nonparametric} questioned its consistency. The interested reader is for instance referred to the introduction of~\cite{derennes2018estimation} for a more detailed discussion on these estimation issues.

In the present paper, the estimation of $\delta_i$ is performed by using the method described in~\cite{derennes2018estimation}: it does not rely on any assumption on the model and works in particular for dependent input. It rests on the copula-representation of $\delta_i$ noted in~\cite{wei2014moment}, namely
\begin{equation}\label{delta_copula}
	\delta_i = \dfrac{1}{2} \int_{0 \leq u, v \leq 1} \left \lvert c_i(u,v)-1 \right \rvert \d u \d v,
\end{equation}
where $c_i$ is the density copula of $(X_i, Y)$, i.e., the density of $(F_{X_i}(X_i), F_Y(Y))$. Based on this representation, the approximation proposed in~\cite{derennes2018estimation} uses a maximum entropy estimation $\hat c_i$ of $c_i$ imposing estimated fractional moments as constraints, and then a Monte Carlo estimation $\frac{1}{2N'} \sum_{k=1}^{N'} \lvert \hat c_i(U^k_1, U^k_2) - 1 \rvert$ of the integral with the $(U^k_1, U^k_2)$ being i.i.d.\ random variables uniformly distributed on $[0,1]^2$.

At this point we stress an important point: all these estimation techniques assume that one can sample from the input distribution $\X$.\ As explained in the introduction however, estimating $\delta^f_i$ amounts to applying these techniques when the input distribution is that of $\X$ conditioned on failure, which is in general unknown. Thus, before applying these methods one needs to be able to sample from $\tilde \X$.

\subsection{Generating conditioned samples {\textnormal{$\tilde{\X}$}}} \label{sub:ASMC}
\subsubsection{General aspects}
The most naive method for generating failure samples is the rejection method. For a given sample $(\X^1, \ldots, \X^N)$ i.i.d.\ with common distribution $f_\X$, a subsample is obtained by recording samples which satisfy $\mathcal{M}(\X^k) > S$. However, this approach leads to a huge computational cost when the failure probability is low.\\ When some information is known on the failure event, this cost can be reduced by leveraging ``good'' auxiliary distributions in importance sampling techniques \cite{bucklew2013introduction}. In reliability, a method widely used for designing auxiliary distributions is shifting the input distribution to a design point, which may be determined thanks to FORM/SORM methods \cite{melchers1990radial}. Importance sampling is then combined with Monte Carlo Markov Chain to generate samples distributed as $\tilde{\X}$ \cite{au2004probabilistic}.\\
Another efficient method to generate conditioned samples $\tilde{\X}$ is the adaptive Sequential Monte Carlo (SMC) procedure proposed and studied in~\cite{cerou2012sequential} that we present in the next section. Several variants have been proposed in different scientific
communities. It was adapted in \cite{au2001estimation} (called subset simulation) for rare event assessment purpose and studied theoretically from the Markov processes point of view in \cite{cerou2012sequential}.\\
 As a final remark, one can mention that importance sampling-based methods and subset simulation may be combined with a surrogate model such as Kriging as it is a powerful tool in the context of costly-to-evaluate computer models. For instance, we can mention the method AK-IS \cite{echard2013combined} which combines Kriging and importance sampling or AK-SS \cite{huang2016assessing} which associates Kriging and subset simulation. However, the counterpart (for the purpose of the present paper) remains the difficulty to catch and measure the impact of the modeling errors induced by the surrogate model itself. \\
We favor in this article the use of the SMC procedure but any of the above mentioned techniques could be applied to generate samples with the same distribution as $\tilde{\X}$. 
\subsubsection{Sequential Monte Carlo}
In what follows, by duplicating a finite set $\{x_k\}$ into $N$, we mean drawing $N$ times independently and uniformly from $\{x_k\}$. The algorithm parameters are $N_x$, $\rho$, $A_x$ and $T$, corresponding respectively to the number of particles, the threshold for the quantile, the number of steps of the Metropolis--Hastings sampler, and the exploration (or proposal) kernel in this sampler.
\begin{description}
	\item[Initialization:] set $p = 0$, generate $(\X^1_p, \ldots, \X^{N_x}_p)$ i.i.d.\ according to $f_{\X}$ and compute $Y^k_p = \mathcal{M}(\X^k_p)$ for $k = 1, \ldots, N_x$;
	\item[Selection:] let $\gamma_p$ be the $\rho$-quantile of the $Y^k_p$: if $\gamma_p > S$, then stop, otherwise duplicate the $\rho N_x$ particles with $Y^k > \gamma_p$ into $N_x$ particles.
	\item[Mutation:] apply $A_x$ times the Metropolis--Hastings algorithm with exploration kernel $T$ and target distribution $\X \mid \mathcal{M}(\X) > \gamma_p$ to each of the $N_x$ particles, denote by $(\X^1_{p+1}, \ldots, \X^{N_x}_{p+1})$ the newly obtained particles with corresponding $Y^k_{p+1} = \mathcal{M}(\X^k_{p+1})$, increment $p$ and go back to the selection step.
\end{description}
The black box is called for every particle at every step of the Metropolis--Hastings sampler in order to compute the acceptance probability, so that if $m$ denotes the (random) number of steps of this algorithm, then the number of calls to the black box $\mathcal{M}$ is equal to $N_x(1 + m A_x)$.

As noted in~\cite{cerou2012sequential}, at the end of this algorithm the $(\X^1_m, \ldots, \X^{N_x}_m)$ are approximately distributed according to $\X \mid Y > \gamma_m$ but are not independent. To improve independence and tune the final size of the sample, an additional step is considered. There are thus two additional parameters, the size $N$ of the sample and the number of steps $A$ of the Metropolis--Hastings sampler in this additional step.
\begin{description}
	\item[Sampling:] duplicate the $N_x$ particles $(\X^1_m, \ldots, \X^{N_x}_m)$ into $N$ particles, and apply $A$ times to each particle the Metropolis--Hastings algorithm with exploration kernel $T$ and target distribution $\X \mid \mathcal{M}(\X) > S$.
\end{description}
This adds $N \times A$ calls to the black box, and the output of this algorithm is a sample $(\tilde \X^1, \ldots, \tilde \X^N)$ which is approximately i.i.d.\ according to $\tilde \X = \X \mid \mathcal{M}(\X) > S$ together with the corresponding values $\tilde Y^k = \mathcal{M}(\tilde \X^k)$.

\subsection{Simultaneous estimation of $\delta^f_i$ and $\bar \eta_i$} \label{sub:simultaneous}

We now explain how to combine the method for estimating $\delta_i$ with the adaptive SMC sampler described above to have a simultaneous estimation of $\delta^f_i$ and $\bar \eta_i$.

\begin{description}
	\item[Step 1 - Input realizations generation.] Using the adaptive SMC procedure of Section~\ref{sub:ASMC}, obtain $(\tilde{\X}^1,\ldots,\tilde{\X}^N)$ approximately i.i.d.\ from $f_{\tilde{\X}}$ and their corresponding value $\tilde Y^k = \mathcal{M}(\tilde \X^k)$ by $\mathcal{M}$.
	\item[Step 2 - Density estimation.] Use the sample $((\tilde X^k_i, \tilde Y^k), k = 1, \ldots, N)$ to obtain estimates $\hat f_{\tilde X_i}$ and $\hat c_i$ of the density $f_{\tilde X_i}$ of $\tilde X_i$ and of the copula $c_i$ of $(\tilde X_i, \tilde Y)$, respectively. In this article, they are both estimated with the maximum entropy method with estimated fractional moments (see~\ref{ME principle}) but any other efficient density and copula estimation technique can be chosen.
	\item[Step 3 - Indices estimation.] Use the estimates $\hat f_{\tilde X_i}$ and $\hat c_i$ to obtain estimates of $\bar \eta_i$ and $\delta^f_i$ as follows:
	\begin{itemize}
		\item for $\bar \eta_i$, estimate the one-dimensional integral $\lVert f_{X_i} - \hat{f}_{\tilde{X}_i} \rVert_{L^1(\R)}$ either by direct numerical approximation, or if $f_{X_i}$ can be sampled from, by Monte Carlo method via
		\[ \hat {\bar \eta}_i = \frac{1}{N'} \sum_{k=1}^{N'} \left \lvert \frac{\hat f_{\tilde X_i}(X^k_i)}{f_{X_i}(X^k_i)} - 1 \right \rvert \]
		where the $X^k_i$ are i.i.d.\ with common distribution $X_i$;
		\item for $\delta^f_i$, generate $((U_1^k,U_2^k), k=1, \ldots, N')$ i.i.d.\ uniformly distributed on $[0,1]^2$ and estimate $\delta_i^f$ by
		\begin{equation} \label{estim_delta_f}
			\Hat{\delta}_i^f = \dfrac{1}{2 N'}\sum_{k=1}^{N'}|\Hat{c_i}(U_1^k,U_2^k)-1|~.
		\end{equation}
	\end{itemize}
\end{description}

It has to be pointed out that the proposed procedure can be applied to output models with correlated inputs even if the interpretation of the results could remain difficult. A potential perspective is to combine this proposed work with Shapley effect estimation of Borgonovo's indices \cite{sarazin2020estimation} as Shapley effects are easier to interpret. \\
As promised, the proposed algorithm also provides simultaneous estimation of both $\delta_i^f$ and $\bar{\eta}_i$ from one common SMC procedure: indeed, after the first step no more call to the black box $\mathcal{M}$ is needed. In particular, the (random) number of calls to the black box is $N_x+m A_x N_x + A N$ as explained in Section~\ref{sub:ASMC} ($N_x+m A_x N_x$ calls for the failure probability estimation and $A N$ calls for sensitivity analysis). If this number of calls is too expensive for a given application, all this procedure can be combined with active learning of a surrogate model for probability estimation such as in \cite{huang2016assessing}, \cite{echard2013combined}. In that case, the $A N$ calls to the black box for sensitivity analysis become $A N$ calls to the surrogate model.

\section{Numerical applications \label{appl}}

In this section, the proposed estimation scheme is applied on four output models. Firstly, we consider two analytical cases for which the unconditional and conditional output distributions are known so that theoretical values of the importance measures $\delta_i^f$ and $\bar{\eta}_i$ are available by using numerical integration. We then consider a single degree of freedom oscillator with $d=6$ independent and lognormally distributed inputs. Finally, as a last test case, we study a launcher stage fallout model which takes $d=6$ input parameters into account.

For each example, computation time and number of model calls are given to assess the efficiency of the proposed method. Results are obtained with a computer equipped with a 3.5 GHz Intel Xeon 4 CPU. 

When the input $\X \sim N(\bm{\nu},\bm{\Sigma})$ is normally distributed, mutation steps in the adaptive SMC algorithm are performed by using the natural exploration kernel so-called Crank Nicholson shaker and defined by 
\begin{equation*}
\label{Crank}
T(\x,\cdot) \sim \textbf{L}\left(\sqrt{1-a}\times\textbf{L}^{-1}(\x-\bm{\nu}) + \sqrt{a}Z\right)+\bm{\nu}~,
\end{equation*}
where $a \in (0,1)$ is a parameter of the kernel, $Z \sim N(0,I_d)$ and $\textbf{L}$ is the lower triangular matrix in the Cholesky decomposition of $\bm{\Sigma}$, i.e., $\bm{\Sigma}=\textbf{L}\textbf{L}^T$.

Standard deviations (Sd) of estimators are computed by performing 100 runs of the proposed scheme in order to study its variability. When a theoretical value $\theta$ is available, the accuracy of an estimator $\hat{\theta}$ is measured by the mean of the relative difference (RD) $\frac{\theta-\hat{\theta}}{\theta}$. 

\subsection{Comparison with ROSA Sobol indices}
It consists in analyzing the influence of the inputs $X_i$ on the variance of $\mathds{1}_{Y>S}$ \cite{luyi2012moment}. Sobol indices on this indicator function can thus be defined in the following way:
$$
S_{i}^{\mathds{1}_{Y>S}} = \frac{\mathbb{V}[{\mathbb{E}[\mathds{1}_{Y>S}}(\mathbf{X})|X_i]]}{\mathbb{V}[{\mathds{1}_{Y>S}(\mathbf{X})}]},
$$
where $\mathbb{V}$ is the variance and $S_{i}^{\mathds{1}_{Y>S}}$ is the first-order Sobol index associated to the variable $X_i$. Advanced sampling-based estimation schemes for these Sobol indices have been investigated in \cite{wei2012efficient} but still require a lot of simulations to achieve convergence. Another efficient estimation procedure using SMC has been proposed recently by \cite{perrin2019efficient} as the first-order $S_{i}^{\mathds{1}_{Y>S}}$ can be rewritten as follows:
$$
S_{i}^{\mathds{1}_{Y>S}}= \frac{\p (Y>S)}{1 - \p (Y>S)} \mathbb{V}\left[{\frac{f_{\tilde X_i}(X_i)}{f_{X_i}(X_i)}}\right].
$$
SMC enables to generate a set of samples $(\tilde{\X}^1,\ldots,\tilde{\X}^N)$ approximately i.i.d.\ from $f_{\tilde{\X}}$ from which it is possible to estimate  a density $\hat f_{\tilde X_i}$ and finally, the Sobol indices are computed with the empirical variance. The following steps describe the complete procedure. 
\begin{description}
	\item[Step 1 - Input realizations generation.] Using the adaptive SMC procedure of Section~\ref{sub:ASMC}, obtain $(\tilde{\X}^1,\ldots,\tilde{\X}^N)$ approximately i.i.d.\ from $f_{\tilde{\X}}$ and estimate $\p (Y>S)$ with $\hat{P}_f$. 
	\item[Step 2 - Density estimation.] Use the sample $(\tilde X^k_i, k = 1, \ldots, N)$ to obtain an estimate $\hat f_{\tilde X_i}$ of the density $f_{\tilde X_i}$ of $\tilde X_i$; In this article, this density is estimated with the maximum entropy method with estimated fractional moments (see~\ref{ME principle}) but any other efficient density estimation technique can be chosen;
	\item[Step 3 - Indice estimation.] Estimate the Sobol indices $S_{i}^{\mathds{1}_{Y>S}}$ from $\hat f_{\tilde X_i}$ and $\hat{P}_f$ as follows:
		\[ \hat S_{i}^{\mathds{1}_{Y>S}} = \left(\frac{\hat{P}_f}{1-\hat{P}_f}\right)\left(\frac{1}{N'} \sum_{k=1}^{N'}  \left(\frac{\hat f_{\tilde X_i}(X^k_i)}{f_{X_i}(X^k_i)}\right)^2 - \left(\frac{1}{N'} \sum_{k=1}^{N'}  \frac{\hat f_{\tilde X_i}(X^k_i)}{f_{X_i}(X^k_i)}\right)^2 \right)  \]
		where the $X^k_i$ are i.i.d.\ with common distribution $X_i$.
\end{description}
There are mainly two limits to this approach: the inputs have to be independent for a correct interpretation of the results and the variance should be a correct indicator of the output variability. Total Sobol indices can also be defined but are more complicated to estimate as they require the estimation of multivariate densities.\\
The estimation of $S_{i}^{\mathds{1}_{Y>S}}$, $\delta^f_i$ and $\bar \eta_i$  is based on the same sample $(\tilde{\X}^1,\ldots,\tilde{\X}^N)$ and thus the computational cost is similar. 

\subsection{Example 1: back to the toy model of the introduction}

We go back to the toy model~\eqref{toy_model} of the introduction, i.e., $Y = X_1 + \mathds{1}_{X_1>S} \lvert X_2 \rvert$ where $S=3$, $X_1$ and $X_2$ are independent, $X_1 \sim N(0,1)$ and $X_2 \sim N(0,5)$. The failure probability can also be evaluated to $1.35 \times 10^{-3}$. We compare in Table~\ref{delta_interpretation} theoretical values with estimates obtained with the proposed method. In average, runs last 317 seconds and make $34{,}640$ calls to the black box ($19{,}460$ calls for the failure probability estimation and $15{,}000$ calls for the sensitivity analysis). From the different relative differences, one can see that $\delta_i^f$ and $\bar{\eta}_i$ estimates are close to their respective reference values and present reasonable variability with regard to the budget allocated to the estimation. Sobol indices $S_{i}^{\mathds{1}_{Y>S}} $ and target sensitivity index $\bar{\eta}_i$ gives very similar results and share thus the same ROSA interpretation. 

\begin{table}[t]
\caption{Estimates of $\delta^f_i$, $\bar \eta_i$ and $S_{i}^{\mathds{1}_{Y>S}} $ of example 1. Set of parameters for the adaptive SMC algorithm: $N_x=500$, $A_x=3$, $\rho = 0.3935$, $a=0.5$, $A = 5$ and $N=3{,}000$.\label{delta_interpretation}}
\centering%
\begin{tabular}{c| c |c c c}
\toprule
\bf{Input} & \bf{Theoretical} & \multicolumn{3}{c}{\bf Estimation $\Hat{\delta}_i^f$} \\
 &\bf{value $\delta_i^f$ (rank)} & Mean (rank) & Sd & RD \\
\midrule
$X_1$ & \bf{0.0781 (2)} & 0.0930 (2) & 0.0101 & -0.1908 \\

$X_2$ & \bf{0.7686 (1)} & 0.7200 (1) & 0.0077 & 0.0632 \\
\midrule
\midrule
\bf{Input} & \bf{Theoretical} & \multicolumn{3}{c}{\bf Estimation $\Hat{\bar \eta}_i$} \\
 &\bf{value $\bar{\eta}_i$ (rank)} & Mean (rank) & Sd & RD \\
\midrule
$X_1$& \bf{0.9987 (1)} & 0.9997 (1) & 0.0095 & -0.001\\
$X_2$ & \bf{0 (2)} & 0.0315 (2) & 0.0103 & \textbf{-}\\
\midrule
\midrule
\bf{Input} & \bf{Theoretical} & \multicolumn{3}{c}{\bf Estimation $\hat S_{i}^{\mathds{1}_{Y>S}} $} \\
 &\bf{value $S_{i}^{\mathds{1}_{Y>S}} $ (rank)} & Mean (rank) & Sd & RD \\
\midrule
$X_1$& \bf{1~(1)} & 1.0225 (1) & 0.0672 & -0.0225\\
$X_2$ & \bf{0 (2)} & ${1.26 \times 10^{-5}~(2)}$ (2) & ${5.53 \times 10^{-6}~(2)}$  & \textbf{-}\\
\bottomrule
\end{tabular}
\end{table} 

\subsection{Example 2: an analytical test case \label{toy_case}}

Let us consider the following output model:
$$Y=X_1+X_2^2$$
where $X_1$ and $X_2$ are i.i.d.\ standard Gaussian random variables and the failure event is $\{Y>15\}$. This model is in the same vein as the previous toy model but slightly more realistic. Unconditional and conditional output distributions are known: $Y \mid X_1$ follows a $\chi^2$-distribution shifted by $X_1$ and $Y \mid X_2$ is normally distributed with unit variance and mean $X_2$. We thus have the following expressions for the densities:
\[ f_Y(y) = \int_0^\infty \dfrac{e^{-\frac{(y-t)^2}{2}-\frac{t}{2}}}{2\pi\sqrt{t}} \d t \]
and
\[ f_{Y\mid X_1}(y) = \frac{e^{-\frac{(y-X_1)}{2}}}{\sqrt{2\pi(y-X_1)}}\mathbb{1}_{y \geq X_1} \ \text{ and } \ f_{Y\mid X_2}(y) = \dfrac{e^{-\frac{(y-X_2)^2}{2}}}{\sqrt{2\pi}} \]
for the conditional densities. Thus, theoretical values of sensitivity measures $(\delta_1,\delta_2)$, $(\delta_1^f,\delta_2^f)$, $(\bar \eta_1, \bar \eta_2)$ and $(S_{1}^{\mathds{1}_{Y>S}}, S_{2}^{\mathds{1}_{Y>S}})$ are available via numerical integration. The failure probability can also be evaluated to $1.2387 \times 10^{-4}$. We gathered in Table \ref{delta_1} the estimates of $\delta_i^f$ and $\bar{\eta}_i$ obtained from the proposed method. In average, runs need 350 seconds to compute all the $\delta$ and $\eta$-indices and make $25{,}200$ calls to the black box ($10{,}200$ calls for the failure probability estimation and $15{,}000$ calls for the sensitivity analysis).

One can see that estimates $\{\hat{\delta}_i^f\}$ respect the good importance ranking, namely $X_2>X_1$. However, the estimation of $\delta^f_1$ exhibits an important difference between average values and reference ones. This difference is due to the fact that the samples $\{\X^k\}$ obtained with the SMC procedure are not completely independent and distributed from $f_{\tilde{\X}}$ since only $A=3$ steps of the Metropolis--Hastings sampler are performed in the final sampling step. 
Indeed, increasing $A$ from $3$ to $30$ leads to average values of $\hat \delta^f_1$ of $0.0206$ with a standard deviation of $0.0091$. 

In this example, the indices $\{\delta_i^f\}$ enable to detect a drastic change in the importance ranking. Indeed, the contribution of the first input $X_1$ becomes negligible at the failure of the system whereas it is the most influential under nominal operation. The indices $\{\bar{\eta}_i\}$ and Sobol indices $\hat S_{i}^{\mathds{1}_{Y>S}} $ lead to the same conclusion, namely that the influence of the input $X_2$ at the failure predominates with $\bar{\eta}_2$ close to $1$.


\begin{table}[t]
\caption{Estimates of $\delta^f_i$, $\bar \eta_i$ and $S_{i}^{\mathds{1}_{Y>S}} $ of example 2. Set of parameters for the adaptive SMC algorithm: $N_x=300$, $A_x=3$, $\rho = 0.5507$, $a=0.5$, $A = 3$ and $N=5{,}000$. \label{delta_1}}
\centering%
\setlength{\tabcolsep}{3pt}
\begin{tabular}{c|c||c|c c c}
\toprule
\bf{Input} & \bf{Theoretical} & \bf{Theoretical} & \multicolumn{3}{c}{\bf Estimation $\Hat{\delta}_i^f$} \\
 &\bf{value $\delta_i$ (rank)} &\bf{value $\delta_i^f$ (rank)} & Mean (rank) & Sd & RD \\
\midrule
$X_1$ & \bf{0.4930 (1)} & \bf{0.001 (2)} & 0.0721 (2) & 0.0266 & -71.1\\

$X_2$ & \bf{0.3049 (2)} & \bf{0.4136 (1)} & 0.3998 (1) & 0.0343 & 0.0334\\
\midrule
\midrule
\bf{Input} & $\backslash$ & \bf{Theoretical} & \multicolumn{3}{c}{\bf Estimation $\Hat{\bar \eta}_i$} \\
& $\backslash$ & \bf{value $\bar{\eta}_i$ (rank)} & Mean (rank) & Sd & RD\\
\midrule
$X_1$ & $\backslash$& \bf{0.2093 (2)} & 0.2066 (1) & 0.0605 & 0.0129 \\
$X_2$ & $\backslash$& \bf{0.9969 (1)} & 0.9723 (2) & 0.0567 & 0.0247\\
\midrule
\midrule
\bf{Input} & $\backslash$ & \bf{Theoretical} & \multicolumn{3}{c}{\bf Estimation $\hat S_{i}^{\mathds{1}_{Y>S}} $} \\
& $\backslash$ & \bf{value $S_{i}^{\mathds{1}_{Y>S}}$ (rank)} & Mean (rank) & Sd & RD\\
\midrule
$X_1$ & $\backslash$& $\bf{4.05 \times 10^{-5}~(2)}$ & ${5.50 \times 10^{-5}~(2)}$ & ${3.35 \times 10^{-5}}$ & -0.3481 \\
$X_2$ & $\backslash$& \bf{0.7074 (1)} & 0.7814 (2) & 0.1012 & -0.1046\\
\bottomrule
\end{tabular}
\end{table} 

\subsection{Example 3: a single Degree of Freedom (SDOF) oscillator} 

In this subsection, a non linear SDOF oscillator \cite{bucher1989time} made of a mass $m$ and two springs with free length $r$ and respective stiffness $c_1$ and $c_2$ is considered. It is subjected to a rectangular load pulse with random duration $t$ and amplitude $F$. The model output is defined as
\[ Y = -3r + \left \lvert \dfrac{2F}{c_1+c_2}\text{sin}\left(\sqrt{\dfrac{c_1+c_2}{m}}\dfrac{t}{2}\right)\right \rvert, \]
i.e., the difference between the maximum displacement response of the system and $3r$. The six input variables $c_1$, $c_2$, $r$, $m$, $t$ and $F$ are assumed to be independent and lognormally distributed with respective parameters given in Table~\ref{distrib_osc}. The failure of the system is achieved when the output $Y$ exceeds the threshold $0$ and the associated failure probability is approximately equal to $9\times 10^{-5}$.

We gathered in Table \ref{indices_1} the estimates of $\delta_i^f$ and $\bar{\eta}_i$ obtained from the proposed method. The $\delta_i$'s are obtained with the method described in \cite{derennes2018estimation}. In average, runs last 960 seconds and make $51{,}725$ calls to the black box ($21{,}725$ calls for the failure probability estimation and $30{,}000$ calls for the sensitivity analysis). It appears that the global importance ranking $X_4 < X_2 < X_5 < X_1 < X_6 < X_3$ drastically differs from the importance ranking provided by the conditional sensitivity indices $\delta^f_i$. Especially, the most influential input $X_3=r$ becomes negligible conditionally on the failure event. Changes are more nuanced as far as target indices are concerned. Indeed, target sensitivity indices $\bar \eta_i$ give approximately the same ranking, except that $X_1$ and $X_2$ predominate.  The Sobol indices $\hat S_{i}^{\mathds{1}_{Y>S}} $ are very low for all the inputs even if they give the same ranking as $\bar \eta_i$. It means that none of the inputs is able alone to reach the failure domain. The estimation of total Sobol indices could be of interest.

\begin{table}[t!]
\begin{center}
\renewcommand{\arraystretch}{0.9}
\caption{Distribution parameters (the mean and the standard deviation of the associated normal distribution) of input variables of the SDOF oscillator.\label{distrib_osc}}
\begin{tabular}{c c c}
\toprule
Input & Mean & Sd \\
\midrule
$c_1$ & 2 & 0.2\\
$c_2$ & 0.2 & 0.02\\
$r$ & 0.6 & 0.05\\
$m$ & 1 & 0.05\\
$t$ & 1 & 0.2\\
$F$ & 1 & 0.2\\
\bottomrule
\end{tabular}
\end{center}
\end{table}
\begin{table}[t!]
\caption{Estimates of $\delta^f_i$, $\bar \eta_i$ and $S_{i}^{\mathds{1}_{Y>S}} $ for the SDOF oscillator. Set of parameters for the adaptive SMC algorithm: $N_x=500$, $A_x=3$, $\rho = 0.4866$, $A = 10$ and $N=3,000$.\label{indices_1}}
\centering%
\begin{tabular}{c|c c | c c | c c | c c }
\toprule
\bf{Input} & \multicolumn{2}{c|}{\bf Estimation $\Hat{\delta}_i$} & \multicolumn{2}{c|}{\bf Estimation $\Hat{\delta}_i^f$} & \multicolumn{2}{c}{\bf Estimation $\Hat{\bar \eta}_i$}& \multicolumn{2}{c}{\bf Estimation $\hat S_{i}^{\mathds{1}_{Y>S}} $}\\
 & Mean (rank) & Sd & Mean (rank) & Sd & Mean (rank) & Sd& Mean (rank) & Sd\\
\midrule
$X_1 = c_1$ & 0.0769 (3) & 0.0066 & 0.0995 (1) & 0.0210 & 0.8332 (2) & 0.0653 & 0.0062 (2) & 0.0062\\

$X_2 = c_2$ & 0.0231 (5) & 0.0050 & 0.0322 (6)& 0.0090 & 0.1352 (5)& 0.0340&  0 (5) & 0\\

$X_3 = r$ & 0.4441 (1) & 0.0063 & 0.0329 (5) & 0.0117 & 0.6494 (3) & 0.0690 & 0.0017 (3) & 0.0014\\

$X_4 = m$ & 0.0219 (6) & 0.0051 & 0.0343 (4) & 0.0101 & 0.1306 (6) & 0.0874 & 0 (6) & 0\\

$X_5 = t$ & 0.0751 (4) & 0.0075 & 0.0474 (3)& 0.0150 & 0.3312 (4)& 0.0710 & 0.0001 (4) & 0\\

$X_6 = F$ & 0.1554 (2) & 0.0074 & 0.0871 (2) & 0.0191 & 0.9078 (1)& 0.0317& 0.0142 (1) & 0.0431 \\
\bottomrule
\end{tabular} 
\end{table} 
\begin{table}[t!]
\caption{Estimates of $\delta^f_i$, $\bar \eta_i$ and $S_{i}^{\mathds{1}_{Y>S}} $ for the SDOF oscillator with a higher budget allocated to the adaptive SMC algorithm. Set of parameters for the adaptive SMC algorithm: $N_x=500$, $A_x=10$ (instead of $3$), $\rho = 0.1813$, $A = 10$ and $N=3,000$.\label{indices_2}}
\centering%
\begin{tabular}{c| c c | c c | c c}
\toprule
\bf{Input} & \multicolumn{2}{c|}{\bf Estimation $\Hat{\delta}_i^f$} & \multicolumn{2}{c}{\bf Estimation $\Hat{\bar \eta}_i$} & \multicolumn{2}{c}{\bf Estimation $\hat S_{i}^{\mathds{1}_{Y>S}} $}\\
 & Mean (rank) & Sd & Mean (rank) & Sd& Mean (rank) & Sd\\
\midrule
$c_1$ & 0.0674 (2) & 0.0150 & 0.7949 (2) & 0.0325& 0.0037 (2) & 0.0017\\

$c_2$ & 0.0275 (5) & 0.0064 & 0.1131 (5) & 0.0173& 0 (5) & 0\\

$r$ & 0.0346 (4) & 0.0089 & 0.5651 (3) & 0.0375& 0.0008 (3) & 0.007\\

$m$ & 0.0267 (6) & 0.0056 & 0.0459 (6) & 0.0196& 0 (6) & 0\\

$t$ & 0.0366 (3) & 0.0074 & 0.2812 (4) & 0.0200& 0.0001 (4) & 0.0001\\

$F$ & 0.1147 (1) & 0.0164 & 0.9205 (1) & 0.0149& 0.0315 (1) & 0.0270\\
\bottomrule
\end{tabular} 
\end{table} 

As in the previous example, variability of obtained estimates is non negligible. 
Here, inputs are lognormally distributed and there is no natural exploration kernel like in the Gaussian case. We can find in \cite{chib1995understanding} a discussion about implementation issues for the choice of the exploration kernel. In the current example, a candidate is drawn by adding a Gaussian noise with the same standard deviation as inputs. With this choice, it appears that we respect standard practice which is to tune the proposal distribution to get around $20\%$--$25\%$ acceptance rate \cite{sherlock2009optimal}. Then, the only way to improve previous results is to increase the budget allocated to Metropolis--Hastings steps by increasing parameter $A$ and decreasing the parameter $\rho$ which regulates values of thresholds involved in the SMC procedure. From Table \ref{indices_2} which displays associated results, one can see that previous observed variability has been reduced. The new computation budget is about $262{,}500$ calls to the model ($232{,}500$ calls for the failure probability estimation and $30{,}000$ calls for the sensitivity analysis), which is quite substantial. Nevertheless, it remains substantially less expensive than the budget required by a classical Monte Carlo procedure. Furthermore, associated computational cost may be reduced by using a surrogate model. For instance, AK-SS method \cite{huang2016assessing} combining Kriging and SMC simulation enables to assess small probabilities while replacing the expensive black box $\mathcal{M}$ by a less time-consuming function.

\subsection{Example 4: a launcher stage fallout model} 
\begin{figure}[htb]
	\begin{center}
		\includegraphics[width=0.80\textwidth]{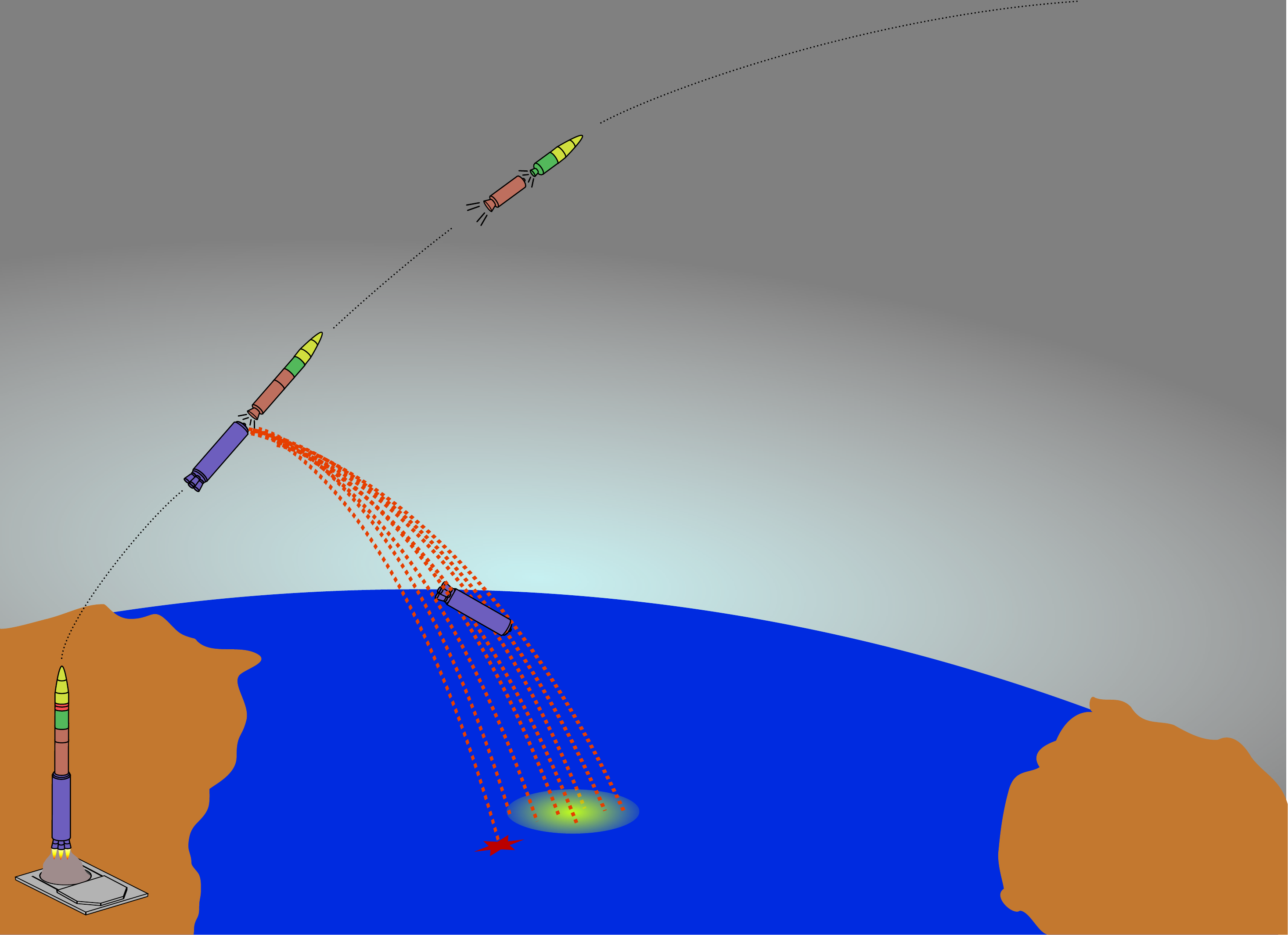}
		\caption{Illustration scheme of a launch vehicle first stage fallout phase into the Atlantic Ocean. Multiple fallout trajectories are drawn (red dotted lines), leading to the safe zone (yellow circular surface). Due to uncertainties, one fallout trajectory may lead to a failure impact point (red star) [cf. \cite{Springer_chapter}].} 
		\label{fig:scenarioV}
	\end{center}
\end{figure}

Space launcher complexity arises from the coupling between several subsystems, such as stages or boosters and other embedded systems. Optimal trajectory assessment is a key discipline since it is one of the cornerstones of the mission success. However, during the real flight, aleatory uncertainties  can affect the different flight phases at different levels (due to weather perturbations, stage combustion etc.) and be combined to lead to a failure state of the space vehicle trajectory. After their propelled phase, the different stages reach successively their separation altitudes and may fall back into the ocean (see Figure~\ref{fig:scenarioV}). Such a dynamic phase is of utmost importance in terms of launcher safety since the consequence of a mistake in the prediction of the fallout zone can be dramatic in terms of human security and environmental impact. That is the reason why it is of prime importance to take it into account during the rare event analysis.

We consider in this section a simplified trajectory simulation code of the dynamic fallout phase of a generic launcher first stage \cite{Springer_chapter} which takes six input parameters into account. The input vector $\X$ contains the following basic variables (i.e., physical variables) representing some initial conditions, environmental variables and launch vehicle characteristics:

\begin{itemize}
\item[] $X_1$: stage altitude perturbation at separation ($\mathrm{\Delta} a$ ($\textrm{m}$));
\item[] $X_2$: velocity perturbation at separation ($\mathrm{\Delta} v$ ($\textrm{m.s}^{-1}$));
\item[] $X_3$: flight path angle perturbation at separation ($\mathrm{\Delta} \gamma$ ($\textrm{rad}$));
\item[] $X_4$: azimuth angle perturbation at separation ($\mathrm{\Delta} \psi$ ($\textrm{rad}$));
\item[] $X_5$: propellant mass residual perturbation at separation ($\mathrm{\Delta} m$ ($\textrm{kg}$));
\item[] $X_6$: drag force error perturbation ($\mathrm{\Delta} C_{d}$ dimensionless).
\end{itemize}

These variables are assumed to be independent and normally distributed with parameters gathered in Table \ref{parameter_lanceur}. As an output, the code will give back the scalar distance $Y = \mathcal{M}(\X)$ which represents the distance between the theoretical fallout position into the ocean and the estimated one due to the uncertainty propagation. The failure event is $\{Y>15\}$ and the associated failure probability is approximately equal to $1.36\times 10^{-4}$.

\begin{table}[h!]
\caption{Input probabilistic model of the launcher phase fallout model.\label{parameter_lanceur}.}
\centering%
\begin{tabular}{@{} c | c c c  @{}}
\toprule
	\bf{Input}  & \bf{Distribution} & \bf{Mean} & \bf{Sd }  \\
	\midrule
	$X_1 = \mathrm{\Delta} a$ ($\textrm{m}$)   & Normal  & $0$  & $165$  \\
	$X_2 = \mathrm{\Delta} v$ ($\textrm{m.s}^{-1}$) & Normal  & 0       & $3.7$  \\
	$X_3 = \mathrm{\Delta} \gamma$ ($\textrm{rad}$)  & Normal  & 0       & $0.001$  \\
	$X_4 = \mathrm{\Delta} \psi$ ($\textrm{rad}$)  & Normal  & $0$            & $0.0018$  \\
	$X_5 = \mathrm{\Delta} m$ ($\textrm{kg}$)  & Normal  & $0$                 & $70$  \\
	$X_6 = \mathrm{\Delta} C_{d}$ ($1$)  & Normal  & $0$                      & $0.1$ \\\bottomrule
\end{tabular}
\end{table}

\begin{table}[h!]
\caption{Estimates of $\delta^f_i$ and $\bar \eta_i$ for the launcher phase fallout model. Set of parameters for the adaptive SMC algorithm: $N_x=800$, $A_x=10$, $\rho = 0.4$, $A = 10$ and $N=5,000$.\label{lanceur_fiab}}
\centering%
\begin{tabular}{c|c c | c c | c c }
\toprule
\bf{Input} & \multicolumn{2}{c|}{\bf Estimation $\Hat{\delta}_i$} & \multicolumn{2}{c|}{\bf Estimation $\Hat{\delta}_i^f$} & \multicolumn{2}{c}{\bf Estimation $\Hat{\bar \eta}_i$}\\
 & Mean (rank) & Sd & Mean (rank) & Sd & Mean (rank) & Sd\\
\midrule
$X_1$ & 0.0156 (5)  & 0.0046 & 0.0149 (5)  & 0.0060  & 0.1235 (4) & 0.0367\\                    

$X_2$ &  0.1535 (2) & 0.0056 & 0.0848 (1) & 0.0095  & 0.8218 (1) &  0.0518\\
 
$X_3$ &   0.0683 (3) & 0.0058 & 0.0406 (2) & 0.0084  & 0.5768 (2) & 0.0514\\                     

$X_4$ &   0.1832 (1) & 0.0050 & 0.0143 (6) & 0.0052  & 0.0722 (6) & 0.0478\\

$X_5$ &  0.0153 (6) & 0.0046 & 0.0160 (4)  & 0.0059  & 0.1230 (5) & 0.0330\\

$X_6$ &  0.0399 (4) & 0.0058 & 0.0284 (3) & 0.0083  & 0.3809 (3) & 0.0660\\
\bottomrule
\end{tabular} 
\end{table} 
Estimates of both target and conditional indices are available in Table \ref{lanceur_fiab}. The $\delta_i$'s are obtained with the method described in \cite{derennes2018estimation}. In average, runs need $750~s$ to compute all the $\delta$ and $\eta$-indices and make $193{,}520$ model calls ($143{,}520$ calls for the failure probability estimation and $50{,}000$ calls for the sensitivity analysis). 

On the one hand, the global ranking provided by $\delta_i$ shows that $X_4$ is the most influential input, followed by $X_2$ and $X_3$. On the other hand, both target and conditional indices underline that the impact of $X_4$ becomes negligible from a rare event perspective. The rest of the ranking remains relatively unchanged and seems to indicate that $X_2$ and $ X_3$ are the most influential input upon the launcher's failure regime. Like the example 2, this test case highlights the importance to perform reliability-oriented sensitivity analysis since the impact of an input parameter (here $X_2$) may depend heavily on the support of the output distribution on which the study is focused. 
    
\section{Generalization: target and conditional sensitivity analysis} \label{sub:generalization}

Following the approach of~\cite{raguet2018target}, we explain here how to generalize our estimation scheme in two directions: $(1)$ considering a more general notion of distance between distributions; $(2)$ assessing the impact of $X_i$ on functions of $Y$.

\subsection{More general distance}

As explained in the introduction, Borgonovo's index is the total variation distance between $(X_i, Y)$ and $(X_i, Y')$ with $Y'$ independent from $X_i$. In the absolutely continuous case, this corresponds to the $L_1$ distance between the joint density $f_{X_i, Y}$ and the product $f_{X_i} f_Y$ of its marginals, which reflects that this index is a measure of dependency between $X_i$ and $Y$. Of course, many other dependency measures exist, for instance the Csisz\'ar dependency measure.

Let $\phi: \R_+ \to \R \cup \{+\infty\}$ be a convex function with $\phi(1) = 0$: then the Csisz\'ar divergence between two probability measures $P$ and $Q$ is given by
\[ \text{div}_\phi \left( P, Q \right) = \int \phi \left( \frac{\d P}{\d Q} \right) \d Q \]
where $P$ is assumed to be absolutely continuous with Radon-Nikodym derivative $\frac{\d P}{\d Q}$ with respect to $Q$. For instance, for $\phi(x) = \frac{1}{2} \lvert 1 - x \rvert$ this is the total variation distance, and for $\phi(x) = - \log(x)$ this is the Kullback--Leibler divergence. From this divergence, we can then define the Csisz\'ar dependency measure (\CDM$_\phi$) between two random variables $Z_1$ and $Z_2$ as
\[ \CDM_\phi(Z_1, Z_2) = \text{div}_\phi \left( (Z_1, Z_2), (Z_1, Z'_2) \right) \]
with $Z'_2$ equal in distribution to $Z_2$ and independent from $Z_1$ (identifying in the above a random variable and its distribution). Because the total variation distance corresponds to the case $\phi(\cdot) = \frac{1}{2} \lvert 1 - \cdot \rvert$, we recover Borgonovo's index with this choice, i.e., we have $\CDM_{\frac{1}{2} \lvert 1 - \cdot \rvert}(X_i, Y) = \delta_i$. Moreover, we note that this dependency measure can still be expressed in a straightforward manner from the copula $c$ of $(Z_1, Z_2)$ provided it exists, namely
\begin{equation} \label{eq:CDM-c}
	\CDM_\phi(Z_1, Z_2) = \int \phi(c(u,v)) \d u \d v,
\end{equation}
thereby generalizing the relation~\eqref{delta_copula} at the heart of our estimation scheme for $\delta^f_i$.

\subsection{Impact on a function of $Y$}

Consider any function $w: \mathcal{M}(\mathds{R}^d) \to \R_+$ such that $w(Y)$ is integrable and let $\tilde \p^w$ be the probability measure which is absolutely continuous with respect to $\p$ with Radon-Nikodym derivative $w(Y)$. Thus, $\tilde \p^w$ is the unique probability measure defined on $(\Omega, \mathcal{F})$ such that
\[ \tilde \p^w(A) = \frac{\E(w(Y) \mathds{1}_A)}{\E(w(Y))} \]
for any measurable set $A \in \mathcal{F}$. 
Adopting the terminology of~\cite{raguet2018target}, we can generalize the two problems laid out in the introduction as follows:
\begin{description}
	\item[Target sensitivity analysis:] what is the influence of $X_i$ on $w(Y)$ (rather than on~$Y$)?
	\item[Conditional sensitivity analysis:] what is the influence of $X_i$ on $Y$ under $\tilde \p^w$ (rather than under $\p$)?
\end{description}

What we have done before corresponds to the case $w(y) = \mathds{1}_{y > S}$. Indeed, for this choice of $w$ the measure $\tilde \p^w \circ \X^{-1}$ is the law of $\tilde \X$ as defined earlier:
\[ \tilde \p^w(\X \in A) = \p(\X \in A \mid Y > S) = \p(\tilde \X \in A). \]

Thus, we generalize $\tilde \X$ to $\tilde \X^w = (X^w_1, \ldots, X^w_d)$ by defining it as a random variable with law $\tilde \p^w \circ \X^{-1}$, and we define $\tilde Y^w = \mathcal{M}(\tilde \X^w)$.

\subsection{Generalization}

In view of the Equations~\eqref{eq:def-eta},~\eqref{eq:def-bar-eta} and~\eqref{eq:def-delta^f} defining $\eta_i$, $\bar \eta_i$ and $\delta^f_i$, respectively, the above extensions suggest the following more general version of these indices:
\[ \eta^{\phi, w}_i = \CDM_\phi \left( X_i, w(Y) \right), \ \bar \eta^{\phi, w}_i = \text{div}_\phi \left( \tilde X^w_i, X_i \right) \ \text{ and } \ \delta^{\phi, w}_i = \CDM_\phi \left( \tilde X^w_i, \tilde Y^w \right). \]

We will assume that $(X_i, Y)$ is absolutely continuous with respect to Lebesgue measure with density $f_{X_i, Y}$, and that $(X_i, w(Y))$ is absolutely continuous with respect to the product measure $\d x \mu(\d a)$ with $\mu$ a measure on $\mathcal{M}(\mathds{R}^d)$ with density $f_{X_i, w(Y)}$. If $w(Y)$ takes values in $\R$, one should typically think of $\mu$ as Lebesgue measure, but this more general formalism also makes it possible to encompass the important case where $w(Y)$ follows a discrete distribution: in this case, $\mu$ should simply be the counting measure and $(X_i, w(Y))$ is automatically absolutely continuous (with respect to $\d x \mu(\d a)$).

Under these assumptions, we have that:
\begin{itemize}
	\item $\eta^{\phi, w}_i = \E \left[ \text{div}_\phi(w(Y), w(Y) \mid X_i) \right]$;
	\item $(\tilde X^w, \tilde Y^w)$ is absolutely continuous with respect to Lebesgue measure with density
	\[ f_{\tilde X^w_i, \tilde Y^w}(x, y) = \frac{w(y) f_{X_i, Y}(x, y)}{\E(w(Y))}. \]
\end{itemize}

For $w(y) = \mathds{1}_{y > s}$ and $\phi(x) = \lvert 1 - x \rvert$, we have the relation~\eqref{eq:eta-bar-eta} between $\eta^{\phi, w}_i$ and $\bar \eta^{\phi, w}_i$ which reads
\[ \eta^{\phi, w}_i = \E(w(Y)) \times \bar \eta^{\phi, w}_i. \]

However, this relation does not seem to hold outside this case, and so in general it is not clear whether $\eta^{\phi, w}_i$ and $\bar \eta^{\phi, w}_i$ can be easily related. Guided by the choice made in the case $w(y) = \mathds{1}_{y > S}$, we consider in the sequel the index $\bar \eta^{\phi, w}_i$ even though it may seem at first glance less natural than $\eta^{\phi, w}_i$.

In order to generalize our estimation scheme, we first need a generalization of the adaptive SMC algorithm of Section~\ref{sub:ASMC}. To sample from the tilted distribution $\tilde \p^w$, usual particle algorithms can be used such as the Metropolis--Hastings sampler with input target density $w(\mathcal{M}(\cdot)) f_\X (\cdot) / \E(w(Y))$. In the case $w(y) = \mathds{1}_{y > S}$ it is hard to sample directly from $\tilde \p^w$ and intermediate distributions, say $\tilde \p^{w_p}$ with $w_p = \mathds{1}_{y > \gamma_p}$, are needed. In this case and with a general $w$, one can for instance use the sequential Monte Carlo samplers proposed in~\cite{del2006sequential}. 

Assume now that one is given a sample $(\tilde \X^{w, 1}, \ldots, \tilde \X^{w, N})$ approximately i.i.d.\ with common distribution $\tilde \X^w$ and their values $\tilde Y^{w, k} = \mathcal{M}(\tilde \X^{w,k})$ by $\mathcal{M}$. As discussed above, in the case $w(y) = \mathds{1}_{y > S}$ this is precisely the purpose of the adaptive SMC algorithm of Section~\ref{sub:ASMC}. Then Step $2$ of our estimation scheme remains unchanged and leads to:
\begin{itemize}
	\item an estimate $\hat f_{\tilde X^w_i}$ of the density $f_{\tilde X^w_i}$ of $\tilde X^w_i$;
	\item an estimate $\hat c^w$ of the copula $c^w$ of $(\tilde X^w_i, \tilde Y^w)$.
\end{itemize}

Using~\eqref{eq:CDM-c} we then have the following two estimations of $\bar \eta^{\phi, w}_i$ and $\delta^{\phi, w}_i$: for $\bar \eta^{\phi, w}_i$, an estimation $\hat {\bar \eta}^{\phi, w}_i$ can be obtained by numerically integrating the one-dimensional integral
\[ \text{div}_\phi(\tilde X^w_i, X_i) = \int \phi \left( \frac{f_{\tilde X^w_i}(x)}{f_{X_i}(x)} \right) f_{X_i}(x) \d x \]
or by a Monte Carlo approximation:
\[ \hat {\bar \eta}^{\phi, w}_i = \frac{1}{N'} \sum_{k=1}^{N'} \phi \left( \frac{f_{\tilde X^w_i}(X^k_i)}{f_{X_i}(X^k_i)} \right) \]
with the $X^k_i$ i.i.d.\ with common distribution $f_{X_i}$. For $\delta^{\phi, w}_i$, draw i.i.d.\ random variables $(U^k_1, U^k_2)$ uniformly distributed on $[0,1]^2$ and consider
\[ \hat \delta^{\phi, w}_i = \frac{1}{N'} \sum_{k=1}^{N'} \phi \left( \hat c^w(U^k_1, U^k_2) \right). \]
\appendix

\section{Maximum entropy principle \label{ME principle}}

\subsection{General principle}

The maximum entropy principle was introduced by Jaynes~\cite{jaynes1957_information}, and the reader is for instance referred to~\cite{kapurentropy} for more details. Let $\mathcal{P}_d(S)$ be the set of probability density functions on $S\subset\mathds{R}^d$, and for $f \in \mathcal{P}_d(S)$ let $H(f)$ be its differential entropy, defined as
\[ H(f) = - \int_{S} \log f(x) f(x) \d x \in [-\infty, +\infty]. \]

In order to choose a density satisfying some constraints $\mathcal{C} \subset \mathcal{P}_d(S)$ (for instance, prescribed first and second moments), the maximum entropy principle asserts to choose among these densities the one with highest entropy, i.e.,
\begin{equation} \label{optim_pb}
	\begin{aligned}
		& \ \underset{f \in \mathcal{P}_d(S)}{\arg \min} & & H(f) \\
		& \ \text{subject to} & & f \in \mathcal{C}
	\end{aligned}
\end{equation}

When the constraints are linear equality constraints, i.e., are of the form $\mathcal{C} = \{f \in \mathcal{P}_d(S): \int \varphi(x) f(x) \d x = \mu\}$ for some $\varphi: \mathds{R}^d \to \mathds{R}^d$ and $\mu \in \mathds{R}^d$, then the above optimization problem is convex and a solution is of the form $f(x) = c e^{-\langle \Lambda^*, \varphi(x) \rangle}\mathds{1}_S(x)$ where $\langle \cdot, \cdot \rangle$ denotes the inner product in $\mathds{R}^d$, $c$ is the normalization constant and $\Lambda^*$ is a feasible solution of the dual optimization problem
\begin{equation} \label{eq:dual-ME}
	\Lambda^* = \underset{\Lambda \in \R^n}{\arg \min} \ \left\{ \langle \Lambda, \mu \rangle + \log \left( \int_S e^{-\langle \Lambda, \varphi(x) \rangle}\d x\right) \right\},
\end{equation}
see for instance~\cite{boyd2004convex} for more details. The above objective function is strictly convex on the set of feasible points and so admits respectively a unique minimum which can been found using standard convex optimization techniques, for instance interior-point algorithms.

The above method can be used to estimate a given density $f_0$: if one knows some moments of the sought density $f_0$, then the idea is simply to put this information as constraints in~\eqref{optim_pb}.

\subsection{Application to Step 2 of our estimation scheme}

In our case, we want to apply the above maximum entropy principle in Step 2 of our estimation scheme (see Section~\ref{sub:simultaneous}) to estimate the density $f_{\tilde X_i}$ of $\tilde X_i$, and the density $c_i$ of $(F_{\tilde X_i}(\tilde X_i), F_{\tilde Y}(\tilde Y))$. Ideally, we would like to consider solutions to~\eqref{optim_pb} with linear equality constraints but the problem is that moments of the sought distributions are unknown. To circumvent this difficulty, we use the sample $((\tilde X^k_i, \tilde Y^k), k = 1, \ldots, N)$ provided by the first step to estimate these moments. Also, for reasons discussed in~\cite{derennes2018estimation} we consider fractional moments for the constraints.

More precisely, consider $\tilde n, n \in \mathds{N}$ and real numbers $\alpha_1 < \cdots < \alpha_{\tilde{n}}$ and $\beta_1 < \cdots < \beta_{n}$, and let
\[ \Hat{M}_{r,s} := \dfrac{1}{N}\sum_{k=1}^N\left(\Hat{F}_{\tilde{X}_i}(\tilde{X}_i^k)\right)^{\alpha_r}\left(\Hat{F}_{\tilde{Y}}(\tilde{Y}^k)\right)^{\alpha_s}, \ r,s=1,\ldots,\tilde{n}, \]
where $\tilde F_{\tilde X_i}$ and $\tilde F_{\tilde Y}$ are the empirical cumulative distribution functions of $\tilde X_i$ and $\tilde Y$, respectively, obtained from the sample $((\tilde X^k_i, \tilde Y^k), k = 1, \ldots, N)$, and
\[ \Hat{M}^i_{t} := \dfrac{1}{N}\sum_{k=1}^N(\tilde{X}_i^k)^{\beta_t},  t=1,\ldots,n. \]

Then the estimates $\hat f_{\tilde X_i}$ and $\hat c_i$ of $f_{\tilde X_i}$ and $c_i$, respectively, are given by
\[ \begin{aligned}
	\hat f_{\tilde X_i} = & \ \underset{f \in \mathcal{P}_1(\supp(\tilde X_i))}{\arg \min} & & H(f) \\
	& \ \text{subject to} & & \int_{\supp(\tilde X_i)} x^{\beta_t} f(x) \d x = \hat M^i_t, \ t = 1, \ldots, n,
\end{aligned} \]
and
\[ \begin{aligned}
	\hat c_i = & \ \underset{f \in \mathcal{P}_2([0,1]^2)}{\arg \min} & & H(f) \\
	& \ \text{subject to} & & \int_{[0,1]^2} x^{\alpha_r} y^{\alpha_s} f(x, y) \d x \d y = \hat M_{r,s}, \ r, s = 1, \ldots, \tilde n.
\end{aligned} \]

These solutions are obtained by the method described above. Note that the number of constraints is then $n$ for estimating $f_{\tilde X_i}$ and $\tilde n^2$ for estimating $c_i$. In this article, $n$ and $\tilde n$ are set to $3$. 

\end{document}